\magnification=1200
\hfuzz=10pt

\font\BBB = msbm10
\def\^#1{\if#1i{\accent"5E\i}\else{\accent"5E #1}\fi}
\def\"#1{\if#1i{\accent"7F\i}\else{\accent"7F #1}\fi}
\def\nn{{\bf n}}
\def\gg{{\bf g}}
\def\hh{{\bf h}}

\def\pp{{\bf p}}
\def\mm{{\bf m}}

\def\NN{{\hbox{\BBB N}}}

\def\ui{\underline{i}}
\def\uj{\underline{j}}

\def\tij{T^{[i,j]}}
\def\dij{D_T^{[i,j]}}

\def\rk{{\rm rk \thinspace}}
\def\wt{{\rm wt \thinspace}}
\def\id{{\rm id \thinspace}}
\def\thick{{\rm thick \thinspace}}
\def\ker{{\rm ker \thinspace}}
\def\card{{\rm card \thinspace}}

\def\sl{{ sl \thinspace}}
\def\SL{{ SL \thinspace}}
\def\Lie{{\rm Lie \thinspace}}
\def\ad{{\rm ad \thinspace}}

\def\cMin{{\rm cMin \thinspace}}
\def\wt{{\rm wt \thinspace}}
\def\p2p{\pi^{\prime \prime}}

\def\CC{\hbox{ \BBB C}}

\def\orb{{\cal O}}
\def\dia{{\cal D}}
\def\rin{{\cal R(V)}}
\def\var{{\cal V}}
\def\war{{\cal W}}
\def\iv{{\cal I(V)}}
\def\p2p{\pi^{\prime \prime}}
\newbox\QEDbox
\newbox\nichts\setbox\nichts=\vbox{}\wd\nichts=2mm\ht\nichts=2mm
\setbox\QEDbox=\hbox{\vrule\vbox{\hrule\copy\nichts\hrule}\vrule}
\def\qed{\leavevmode\unskip\hfil\null\nobreak\hfill\copy\QEDbox\medbreak}
\ifx\Tableauxloaded\relax \endinput \else \let\Tableauxloaded\relax \fi

% Dessiner des tableaux de Young

\newcount\cols
{\catcode`,=\active\catcode`|=\active
 \gdef\Young(#1){\hbox{$\vcenter
 {\mathcode`,="8000\mathcode`|="8000
  \def,{\global\advance\cols by 1 &}%
  \def|{\cr
        \multispan{\the\cols}\hrulefill\cr
        &\global\cols=2 }%
  \offinterlineskip\everycr{}\tabskip=0pt
  \dimen0=\ht\strutbox \advance\dimen0 by \dp\strutbox
  \halign
   {\vrule height \ht\strutbox depth \dp\strutbox##
    &&\hbox to \dimen0{\hss$##$\hss}\vrule\cr
    \noalign{\hrule}&\global\cols=2 #1\crcr
    \multispan{\the\cols}\hrulefill\cr%
   }
 }$}}
 \gdef\Skew(#1:#2){\hbox{$\vcenter
 {\mathcode`,="8000\mathcode`|="8000
  \dimen0=\ht\strutbox \advance\dimen0 by \dp\strutbox
  \def\boxbeg{\vbox
    \bgroup\hrule\kern-0.4pt\hbox to\dimen0\bgroup\strut\vrule\hss$}%
  \def\boxend{$\hss\egroup\hrule\egroup}%
  \def,{\boxend\boxbeg}%
  \def|##1:{\boxend\vrule\egroup\nointerlineskip\kern-0.4pt
    \moveright##1\dimen0\hbox\bgroup\boxbeg}%
  \def\\##1\\##2:{\boxend\vrule\egroup\nointerlineskip\kern-0.4pt
    \kern ##1\dimen0\moveright##2\dimen0\hbox\bgroup\boxbeg}%
  \moveright#1\dimen0\hbox\bgroup\boxbeg#2\boxend\vrule\egroup
 }$}}
}

% We need script and scriptscript sizes of text italic in small squares

\font\sevenit=cmti7 \font\fiveit=cmti5
\scriptfont\itfam=\sevenit \scriptscriptfont\itfam=\fiveit

\def\smallsquares
{\textfont0=\scriptfont0 \scriptfont0=\scriptscriptfont0
 \textfont1=\scriptfont1 \scriptfont1=\scriptscriptfont1
 \textfont\itfam=\scriptfont\itfam \scriptfont\itfam=\scriptscriptfont\itfam
 \textfont\bffam=\scriptfont\bffam \scriptfont\bffam=\scriptscriptfont\bffam
 \setbox0=\hbox{$($}
 \setbox\strutbox=\hbox{\vrule width 0pt height\ht0 depth\dp0 }% smaller
}

\def\bigsquares
{\setbox\strutbox=\hbox
 {\vrule width 0pt height1.3\ht\strutbox depth1.3\dp\strutbox}% larger
}

\newbox\hpair \newbox\vpair
\def\installboxes
{\dimen0=\ht\strutbox \advance\dimen0 by \dp\strutbox
 \setbox0=\vbox{\hrule width \dimen0}
 \setbox1=\rlap{\strut\vrule}
 \setbox\hpair=\rlap{\raise\ht\strutbox\copy0 \lower\dp\strutbox\copy0}
 \setbox\vpair=\rlap{\raise\dimen0\copy1 \kern\dimen0 \copy1}
}

\newcount\ribbonlength
\newcount\joff \newcount\jmax \newcount\imid \newcount\jmid
\def\ribbon #1,#2:#3;#4
{\setbox2=\hbox
 {\count0=0 \count1=0
  \copy1 \rlap{\lower\dp\strutbox\copy0}\if|#4|\else\ribs#4 \fi
  \raise\count1\dimen0\rlap
    {\kern\count0\dimen0 \raise\ht\strutbox\copy0 \copy1\raise0.4pt\copy1}%
  \global\joff=\count0
  \divide\dimen0 by 2 \raise\imid\dimen0\rlap{\kern\jmid\dimen0
    \hbox to 2\dimen0{\hfil$#3$\hfil}}%
 }%
 \lower#1\dimen0\rlap{\kern#2\dimen0\box2}%
 \advance\joff by #2\relax \ifnum \joff>\jmax \global\jmax=\joff \fi
}
\def\ribs#1#2
{\raise\count1\dimen0\rlap
   {\kern\count0\dimen0 \copy\ifnum #1=0 \hpair \else \vpair \fi}%
 \advance\ribbonlength by -2
 % Now increment \count0 or \count1, depending on #1. Also, if
 % \ribbonlength=1, set (\jmid,\imid):=2*(\count0,\count1) after
 % incrementation, and if \ribbonlength=0 set (\jmid,\imid) to the sum of
 % (\count0,\count1) before and after incrementation. Since `2\count0' is not
 % a valid <integer>, the easiest way is to set (\jmid,\imid) explicitly both
 % for \ribbonlength=0 and \ribbonlength=1, and to then advance the
 % coordinates if \ribbonlength=0 OR \ribbonlength=1
 \ifnum \ribbonlength=0 \jmid=\count0 \imid=\count1 \fi
 \advance\count#1 by 1
 \ifnum \ribbonlength=1 \jmid=\count0 \imid=\count1 \fi
 \ifnum \ribbonlength>-1 \ifnum \ribbonlength<2
   \advance \jmid by \count0 \advance \imid by \count1
 \fi\fi
 \if |#2|\else \ribs#2 \fi % loop while there are bits in the specified string
}

\def\rtab#1
{\installboxes \ribbonlength=#1 \global\jmax=0
 \setbox0=\hbox{\strut\aftergroup\endrtab \aftergroup}}
\def\endrtab{\advance\dimen0 by \jmax\dimen0 \wd0=\dimen0 \vcenter{\box0}}

\def\arrow#1,#2:#3
{\smash{\lower#1\dimen0\rlap{$\kern#2.5\dimen0
 \if r#3 \,\rightarrow \else
 \if l#3 \llap{$\leftarrow\,$} \else
 \if u#3 \setbox0=\hbox to 0pt{\hss$\uparrow$\hss}\ht0=0pt 
         \raise.5\dimen0 \box0 \else
 \if d#3 \setbox0=\hbox to 0pt{\hss$\downarrow$\hss}\dp0=-2pt
         \lower.5\dimen0 \box0 \else
 \if i#3 \setbox0=\hbox to 0pt{\hss$\nwarrow\,$}\ht0=0pt
         \raise.5\dimen0 \box0 \else
 \if o#3 \setbox0=\hbox to 0pt{$\,\searrow$\hss}\dp0=-2pt
         \lower.5\dimen0 \box0 \else
 \fi\fi\fi\fi\fi\fi
 $}}%
}

\centerline{{\bf On the hypersurface orbital varieties of $\sl(N,{\CC})$}}

\medskip
\centerline{Elise Benlolo\footnote{$^{\dag}$}{D\'epartement de Math\'ematiques, Moulin de la Housse, Universit\'e de Reims, F-51687 Reims, France, email: Elise.Benlolo@univ-reims.fr} and Yasmine B. Sanderson\footnote{$^{\ddag}$}{Department of Mathematics, William Paterson University, Wayne, N.J. 07470 USA, email: sandersony@wpunj.edu}}

\bigskip\noindent
{\bf Abstract:} We study the structure of hypersurface orbital varieties of $\sl(N,{\CC})$ (those that are hypersurfaces in the nilradical of some parabolic subalgebra) and how information about this structure is encoded in the standard Young tableau associated to it by the Robinson-Schensted algorithm. We present a conjecture for the exact form of the unique non-linear defining equations of hypersurface orbital varieties and proofs of the conjecture in certain cases.

\bigskip\noindent
{\bf I. Introduction}

\noindent
Let $G$ be a complex semi-simple algebraic group with
Lie algebra $\gg$, on which it acts through the adjoint
representation. A $G$-orbit $\orb$ in $\gg$ is said to be {\it
  nilpotent} if it consists of 
nilpotent elements. Fix some Cartan decomposition of $\gg = \nn_{-}
\oplus \hh \oplus \nn$.  Then, an irreducible component of $\orb
\cap \nn$ is called an {\it orbital variety}. These varieties figure prominently in the primitive ideal theory of ${\cal U}(\gg)$ and the ongoing attempt to establish an ``orbit correspondence'' for semisimple groups (see [B], [BV], [Mc], [J3] for example). In the case of orbital varieties, it was shown by Spaltenstein [Sp] and Steinberg [St2] that the dimension of an orbital variety is half the dimension of the corresponding nilpotent orbit. Joseph [J2] showed that this implies that orbital varieties are Lagrangian. In the orbit method, one would wish to find a correspondence between these Lagrangian subvarieties of (co-)adjoint orbits and simple highest weight modules. Noting that the unions of the closures of orbital varieties arise as associated varieties of simple highest weight modules, Joseph [J4] laid out a program of ``quantization'': He called an orbital variety $\var$ {\it weakly quantizable} if its closure is the associated variety of a simple highest weight module. An orbital variety $\var$ is {\it strongly quantizable} if there exists a highest weight module $M$ whose formal character (as $\hh$-module) matched that of the coordinate ring of $\overline{\var}$ (see [Be]). Benlolo gave two examples of varieties in $sl(6)$ which were strongly quantizable, but only by non-simple highest weight modules [Be]. Melnikov showed that every variety in $sl(N)$ is weakly quantizable [M3]. In [J4] Joseph studies the orbital varieties in the minimal nonzero orbit for a complex semisimple Lie algebra $\gg$ and shows that every such orbit contains at least one strongly quantizable variety. However, he also finds examples of varieties that are not weakly quantizable and varieties that are weakly but not strongly quantizable.

One stumbling block in the study of orbital varieties and related highest weight modules is simply that the structure of orbital varieties remains quite mysterious. Except in the (obvious) case of the Richardson varieties (whose defining relations are all linear) there are no general formulas for the defining equations of orbital varieties. This, for one, makes studying the character of the coordinate ring of $\overline{\var}$ rather difficult if not impossible. Obtaining an exact description of the ideal of definition of an orbital variety would also greatly benefit the calculation of the {\it characteristic polynomial} $p_{\var} \in S(\hh^*)$ of $\var$.   The importance of characteristic polynomials is
revealed through their many characterizations. As $\var$ runs over
the components of $\orb \cap \nn$, where $\orb$ is a fixed orbit, the
$p_{\var}$ 
span a $W$-submodule of $S(\hh)$. This is the representation of $W$
assigned to $\orb$ by the Springer
correspondence ([J2], [Ho]). In addition, the $p_{\var}$ are intimately connected to Goldie rank polynomials [J2] and can also be viewed as
equivariant characteristic classes of 
orbital cone 
bundles [BBM]. Characteristic polynomials can be calculated from the character of the coordinate ring of $\var$ [J2] or directly from a recursive algorithm [J1]. However, knowing the ideal of definition would greatly help in converting the theory into practice.

Our interest in orbital varieties comes mainly from a combinatorial point of view. We therefore restrict our attention to the orbital varieties of $\sl(N,\CC)$: through the Robinson-Schensted correspondence, the set of orbital varieties so $sl(N, \CC)$ is in bijection with the set of standard Young tableaux with $N$ boxes. This bijection is ``natural'' in the sense that information about an orbital variety $\var$ can be ``read off'' the associated tableau $T$. From a standard Young tableau, one can determine the orbit in which an orbital variety lies, its dimension and its $\tau$-invariant $\tau(\var)$, a certain subset of the set $\Pi$ of simple roots. From $\tau(\var)$ one knows the maximal parabolic subgroup $P_{\tau}$ of $\SL(N)$ which stabilizes $\var$: it is generated by the Borel subgroup $B$ and the root vectors $X_{-\alpha}$ where $\alpha \in \tau(\var)$. 

Since an orbital variety is, in some sense, determined by a standard Young tableau, one would like to be able to obtain more information about the structure of $\var$ directly from the combinatorial information in its associated tableau.
With this idea in mind, we concentrated our efforts on the {\it hypersurface orbital varieties} of $\sl(N,\CC)$, that is, orbital varieties which are hypersurfaces in the nilradical of some parabolic. The equations of these varieties are all linear except for one, $f=0$, where $f$ is a homogeneous polynomial in $S(\nn)$ with $\deg(f) \geq 2$. Clearly the linear equations are all of the form $X_{\alpha} = 0$ where $\alpha$ is a sum of simple roots in $\tau(\var)$. So the real problem was extracting information about $f$ from $T$. 

Our main idea was to compare the tableau $T$ with $T_R$, the standard Young tableau associated to the Richardson orbital variety $\var_R$ with the same $\tau$-invariant as $\var$. The relationship between $\var$ and $\var_R$ is the following: $\var \subset \overline{\var}_R = m_{\tau(\var)}$ where $ m_{\tau(\var)}$ is the nilradical of the parabolic subalgebra $\Lie(P_{\tau(\var)})$. This relationship between varieties translates to the following relationship between tableaux: $T$ is obtained from dropping one box of $T_R$ down one row. This allowed us to determine the minimal connected subset $\sigma$ of $\Pi$ such that $f \in S(\nn_{\sigma})$ where $\nn_{\sigma}$ is the subalgebra of $\nn$ generated by $X_{\alpha}$ with $\alpha \in \sigma$. In other words, we can use $T$ to tell us ``where $f$ is located''. An important tool in our proofs are the so-called {\it power-rank conditions} of van Leeuwen [vanL]. These conditions describe relations among the coordinates of a generic matrix in terms of the shapes of the subtableaux of $T$. We were then able to describe $f$ explicitly for many cases of hypersurface orbital varieties. We also conjecture an explicit formula for $f$ when $\var$ a is an arbitrary hypersurface orbital variety. Using previous results of [J1], [BBM] we obtain an explicit formula for the characteristic polynomial of such varieties. We provide extensive examples to illustrate our points.

\noindent {\bf Acknowledgment:}
We would like to thank A. Joseph for inspiring us to work on orbital varieties and the referee for making useful comments which improved the final version of this paper.

\bigskip \noindent
{\bf II. Some background and notation}

Let $G := \SL(N, \CC)$ and let $\gg = \sl(N, \CC)$ denote its Lie algebra. Let $\gg = \nn_{-} \oplus \hh \oplus \nn$ be the Cartan decomposition where $\hh$ denotes the Cartan subalgebra and $\nn$ the nilpotent subalgebra of strictly upper-triangular matrices. Let $\Pi$ denote the set of $N-1$ simple roots $\alpha_1, \ldots, \alpha_{N-1}$. For $i \leq j$, let $\alpha(i,j)$ denote the positive root $ \alpha_i + \alpha_{i+1} +\alpha_{i+2} + \cdots + \alpha_j$. Let $R^+ := \{ \alpha(i,j) \mid 1 \leq i \leq j \leq N-1\}$, the set of positive roots.  Let $X_{\alpha(i,j)} \in \nn$ denote the associated root vector. Note that a generic matrix $x$ in the one-dimensional space spanned by $X_{\alpha(i,j)}$ satisfies $x_{m,n} = 0$ for $(m,n) \not= (i,j+1)$. We will denote by $x_{i,j+1}$ the coordinate corresponding to the root $\alpha(i,j)$. This allows us to identify $S(\nn) \cong \CC[x_{ij} \mid 1 \leq i < j \leq N]$.

 Let $s_i$ denote the simple reflection with respect to the simple root $\alpha_i$ ( $1 \leq i \leq N-1$). We also use $s_i$ to denote the transposition $(i \ i+1)$. The $\{s_i\}_{i=1}^{N-1}$ generate the Weyl group $W$ which, in the case of $\SL(N)$, is isomorphic to $S_N$, the symmetric group on $N$ letters.
For any subset $\omega \subseteq \Pi$, let $R^+(\omega)$ denote the positive roots in $\Pi$ which are sums of simple roots of $\omega$. Let $W(\omega)$ be the Weyl group 
element generated by the simple reflections $s_{\alpha}$ for $\alpha
\in \omega$.

\bigskip \noindent
{\bf III. Orbital varieties and Young tableaux}

$G = \SL(N, \CC)$ acts on $\gg$ by conjugation. When $X \in \nn$, the orbit $\orb := G\cdot X$ is called a {\it nilpotent orbit}. The set of nilpotent orbits of $\sl(N, \CC)$ is in bijection with the set of partitions $\lambda$ of $N$. In fact, to each partition $\lambda=(\lambda_1,\ldots,\lambda_n)$, one can associate the strictly upper-triangular nilpotent matrix $x_{\lambda}$ with $n$ Jordan blocks of size $\lambda_1, \lambda_2,\ldots, \lambda_n$. Thus, the orbit $\orb_{\lambda} := G\cdot x_{\lambda}$ consists of all nilpotent matrices with Jordan canonical form $x_{\lambda}$.

An {\it orbital variety} is an irreducible component of the intersection $\orb \cap \nn$. A more explicit general description of these varieties, due to [St1], [J2], is as follows: Let $ w \in W$ be a Weyl group element. Let $B$ denote the Borel subgroup of $\SL(N)$. Set $\nn \cap w(\nn) = \oplus_{\alpha \in R^+ \cap w^{-1}(R^+)} \CC X_{\alpha}$. Then 
$$\var = \var(w) := (\overline{B \cdot (\nn \cap w(\nn))}) \cap \orb$$
is an orbital variety and the map $w \mapsto \var(w)$ is a surjection of $W$ onto the set of all orbital varieties of $\sl(N)$.

Since the set of nilpotent orbits is indexed by partitions, it is natural to wonder if this indexing somehow extends to orbital varieties. Such an extension exists, which we now describe. We can identify a partition $\lambda$ with a {\it Young diagram} consisting of $N$ boxes with $\lambda_1$ boxes in the $1^{st}$ row, $\lambda_2$ boxes in the $2^{nd}$ row and so on. A {\it standard Young tableau} is a filling of the $N$ boxes with the numbers $1,2, \ldots, N$ in such a way that the numbers increase from left to right in every row and from top to bottom in every column.

\medskip \noindent
{\bf Example 1:} The Young diagram associated to the partition $(4,2,1)$ is below on the left. On the right are several examples of standard Young tableaux of shape $(4,2,1)$.
$$\Young(,,,|,|) \qquad \qquad \Young(1,3,5,6|2,7|4)\qquad \Young(1,2,3,4|5,6|7)\qquad \Young(1,3,6,7|2,4|5)$$

\medskip 
\noindent
{\bf Theorem 1:} [J2, 9.14] The set of orbital varieties in the nilpotent orbit $\orb_{\lambda}$ is in bijection with the set ${\cal T}_{\lambda}$ of standard Young tableaux of shape $\lambda$.

This bijection is a corollary of the Robinson-Schensted correspondence (see [M1], [vanL] for nice descriptions), which associates to each permutation $w \in W$ a certain pair of standard Young tableaux $(A(w), B(w))$. The tableau which will be associated to the orbital variety $\var(w)$ is the tableau $B(w)$ associated to $w$ by this bijection.

Extensive research has been done studying this connection between orbital varieties. From a standard Young tableau $T$, one
can read information about the associated orbital variety $\var \subset \orb_{\lambda}$. In particular, if $\lambda^{\prime}$ is the dual partition of $\lambda$, then $\dim \var = {1 \over 2}(N^2 - ({\lambda^{\prime}}^2_1 + {\lambda^{\prime}}^2_2 + \cdots + {\lambda^{\prime}}^2_j))$ [SS]. From the tableau $T$, one can also determine the $\tau$-invariant $\tau(\var)$ of an orbital variety $\var$.
By definition, $\tau(\var)\subseteq \Pi$ is the set of all simple roots $\alpha$ such that the subgroup $P_{\alpha} := \langle \exp \ad X_{\alpha},\  \exp \ad X_{-\alpha} \rangle$ stabilizes $\var$.
In other words, if $P$ is the maximal parabolic subgroup which stabilizes $\var$, then $P = P_{\tau} := \langle B,\  \exp \ad X_{-\alpha} \ \mid \  \alpha \in \tau(\var) \rangle$.
It turns out that $i$ is above $i+1$ in $T$ if and only if $\alpha_i \in \tau(\var)$ [Ja].

To every subset $\tau \subseteq \Pi$, there exists a (unique) orbital variety $\var_R$ of maximal dimension whose $\tau$-invariant is $\tau$. $\var_R$ is a {\it Richardson variety}, that is $\dim(\var_R)$ equals the dimension of the nilradical $$\mm_{\tau} := \bigoplus_{\alpha \in R^+ \setminus R^+(\tau)} \CC X_{\alpha}$$ of the parabolic subalgebra $\pp_{\tau} = \Lie(P_{\tau})$. Therefore its standard Young tableau $T_R$ is ``top-heavy''. One constructs it by putting the numbered boxes in the topmost row possible such that the restrictions imposed by $\tau$ are respected. 

\medskip\noindent
{\bf Example 2:} Let $G = \SL(8)$. Then each Young diagram will consist of 8 boxes. The standard Young tableau associated to a Richardson orbital variety $\var_R$ with $\tau$-invariant $ \{ \alpha_2, \alpha_3, \alpha_7 \}$ is
$$T_R = \Young(1,2,5,6,7|3,8|4)$$
This orbital variety lies in the orbit $\orb_{\lambda}$ where $\lambda = (5,2,1)$. Then 
$\lambda^{\prime} = (3,2,1,1,1)$. We have that
$$\dim(\var_R) = {1\over 2}(8^2 - (3^2 + 2^2 + 1^2 + 1^2+ 1^2)) = 24 = 28-4 = \card(R^+) - \card(R^+(\tau))=\dim \mm_{\tau} .$$
Let $T_R$ be the standard Young tableau associated to the Richardson variety $\var_R$. A {\it chain} $C$ of $T_R$ is an invariant subset of $\{1,2,\ldots,N \}$ under the action of $W(\tau(\var_R))$. In other words, it is a set of the form $C = \{i, i+1, \ldots, i+k\}$  $(k \in {\bf N} )$, where 

 (a) $i$ is in the first row of $T_R$ 

 (b) if $i+1$ is also in the first
row of $T_R$, then $C=\{i\}$, i.e. $k=0$ 

 (c) if $i+1$ is not in the first row of $T_R$, one requires that 
$\{\alpha_{i+1},...,\alpha_{i+k-1}\} \subseteq \tau$ and, whenever
$i+k<N$, $\alpha_{i+k} \not\in \tau$. 

We say that $C$ has length $k+1$ and denote this by $\mid C\mid = k +1$. Notice that $T_R$ is completely determined by its chains. If $T_R$ is of shape $\lambda = (\lambda_1,\ldots,\lambda_r)$, then it has $\lambda_1$ chains and the number of columns of length $i$ equals the number of chains of length $i$. 

\medskip
\noindent
{\bf Example 3:} The tableau in our previous example has five chains: $C_1 = \{ 1\}$, $C_2=\{2,3,4\}$, $C_3=\{5\}$, $C_4=\{6\}$, $C_5=\{7,8\}$ which are exactly the invariant subsets of $\{1,2,\ldots,8\}$ under the action of the transpositions $s_2=(2 \enspace 3)$, $s_3 = (3 \enspace 4)$ and $s_7 = (7 \enspace 8)$.

\medskip
Consider the orbital varieties $\var_1$ and $\var_2$. We say that
$\var_2$ is a {\it descendant} of $\var_1$ if $\var_2 \not= \var_1$, $\var_2 \subset \overline{\var}_1$ and if any orbital variety satisfies $\war$ satisfies $\var_2 \subseteq \overline{\war} \subseteq \overline{\var}_1$, then either $\var_2 = \war$ or $\war = \var_1$. If $T_1$ (resp. $T_2$) is the standard Young tableau associated to $\var_1$ (resp. $\var_2$), then we say that $T_2$ is a descendant of $T_1$ if and only if $\var_2$ is a descendant of $\var_1$.  

We now consider a hypersurface orbital variety $\var \subset \overline{\var}_R$ with the same $\tau$-invariant. Then $\var$ is a descendant of $\var_R$. Let $T$ (resp. $T_R$) be the standard Young tableau associated to $\var$ (resp. $\var_R$). For any tableau $T$, let $r_T(i)$ denote the row of $T$ in which the box numbered $i$ is located.

\medskip\noindent
{\bf Lemma 1:} $T$ is obtained from $T_R$ by moving a box containing the maximal element of some chain $C$ of length $i$ from row $i$ to row $i+1$.

\smallskip\noindent
{\bf Proof:}
 Let $T_R \in {\cal T}_{\lambda}$ and $T \in {\cal T}_{\mu}$. It is a result of Gerstenhaber that $\var \subset \overline{\var}_R$ implies that $\mu \leq \lambda$ (see [He]). Recall that $\mu \leq \lambda$ is defined as $\mu_1 + \mu_2 + \cdots \mu_i \leq \lambda_1 + \lambda_2 + \cdots + \lambda_i$ for all $i$. Suppose first that $\mu$ differs from $\lambda$ by the dropping of one box (down one or possibly several rows). If $\lambda^{\prime} = ( \lambda_1^{\prime}, \ldots, \lambda_r^{\prime})$ and $\mu^{\prime} = ( \mu_1^{\prime}, \ldots, \mu_s^{\prime})$, then either $r=s$ or $s=r-1$. In the latter case, we will set $\mu_r^{\prime} = 0$. Then there exist $j$ and $k$ ($j < k$) such that $\mu_j^{\prime} = \lambda_j^{\prime} + 1$, $ \mu_k^{\prime} = \lambda_k^{\prime} - 1$ and $ \mu_i^{\prime} = \lambda_i^{\prime}$ when $i \not= j,k$. So 
$$\eqalign{\dim \var &= {1\over 2}(N^2 - (\mu_1^{\prime})^2 - \cdots - (\mu_r^{\prime})^2) \cr
& = {1\over 2}(N^2 - (\lambda_1^{\prime})^2 - \cdots - (\lambda_j^{\prime})^2- \cdots - (\lambda_k^{\prime})^2- \cdots -(\lambda_r^{\prime})^2) - (\lambda_j^{\prime} - \lambda_k^{\prime} +1)\cr
& = \dim \var_R - (\lambda_j^{\prime} - \lambda_k^{\prime} +1)\cr}$$
which equals $\dim \var_R - 1$ if and only if $\lambda_j^{\prime} = \lambda_k^{\prime}$. This can only happen if the box was knocked down one row. Since $\dim \var = \dim \var_R - 1$, then $\mu$ can not be obtained from $\lambda$ by moving more than one box.

  Notice that, since $\tau(\var) = \tau(\var_R)$, then $r_T(i) \geq r_{T_R}(i)$ for all $i$. Since only one box is moved from $T_R$ in order to obtain $T$, it must contain the maximal element of some chain in $T_R$. (Any ``shuffling'' of the boxes would either produce something that is not a standard tableau or would change the $\tau$-invariant).
 \qed

\medskip\noindent
{\bf Example 4:} Consider the Young tableau from Example 2. By knocking down box number $5$ from the first to the second row, one obtains a Young tableau associated to a hypersurface orbital variety contained in $\overline{\var}_R$:
$$\Young(1,2,5,6,7|3,8|4) \enspace \rightarrow \enspace \Young(1,2,6,7|3,5,8|4)$$
The corresponding orbital variety $\var$ is contained in $\orb_{\mu}$ where $\mu = (4,3,1)$. Since $\mu^{\prime} = (3,2,2,1)$, we have that $\dim(\var) = 23 = \dim(\var_R) - 1$. Hence $\var$ is a hypersurface variety. Notice that there are no other ways that one could move a box down one row without changing the $\tau$-invariant. Therefore, in this case, $\overline{\var}_R$ contains only one hypersurface orbital variety with the same $\tau$-invariant.

\bigskip
Likewise, given a standard Young tableau associated to a hypersurface orbital variety $\var$, one can always obtain the tableau associated to the Richardson variety which contains $\var$ by moving an appropriate box up one row.

\bigskip \noindent
{\bf IV. Subtableaux and projections of orbital varieties}

In the following, for any $i$, we denote by $\nn_i$ the subalgebra of strictly upper-triangular matrices in $sl(i)$.

Let $\pi_{1,N-1}: \nn_N \rightarrow \nn_{N-1}$ be the projection which, to a generic matrix $x \in \nn_N$, assigns the same matrix with the $N^{\rm th}$ row and column removed. Let $\var$ be an orbital variety with standard Young tableau $T$. It results from work of Sch\"utzenberger, Knuth and Melnikov (See [M1] Lemma 1.1.3, Theorems 1.3.13 and 4.1.2) that $\pi_{1,N-1}(\var)$ is dense in $\overline{\war} \subset \nn_{N-1}$ where $\war$ is a certain orbital variety. The standard Young tableau $T^{[1,N-1]}$ associated to $\war$ is obtained from $T$ by removing the box with the largest entry. (See also [vanL] for a discussion of this in terms of flag varieties.)

\medskip
\noindent
{\bf Example 5:} If $\displaystyle{ T = \Young(1,2|3,4)}$ then $\displaystyle{T^{[1,N-1]} =  \Young(1, 2|3) }$ .

\medskip
Likewise, let $\pi_{2,N}: \nn_N \rightarrow \nn_{N-1}$ be the projection which, to a generic matrix $x \in \nn_N$, assigns the same matrix with the $1^{\rm st}$ row and column removed. Let $\var$ be an orbital variety with standard Young tableau $T$. In the same way, we will associate to this projection a certain standard Young tableau $T^{[2,N]}$, obtained in the following way: Apply to $T$ the Sch\"utzenberger ``jeu de taquin'' algorithm (see [M1] or [vanL], \S 4): remove the box in the first row and first column to leave an empty square in its place. Then the following step is repeated until the empty square is a corner of the original tableau: move into the empty square the smaller of the entries located directly to the right of and below it. Replace each of the entries $i$ in this tableau by $i-1$ to obtain a standard Young tableau which we will denote by $T^{[2,N]}$. We have that $\pi_{2,N}(\var)$ is dense in $\overline{\war} \subset \nn_{N-1}$ where $\war$ is the orbital variety with associated standard Young tableau $T^{[2,N]}$.

\smallskip\noindent
{\bf Example 6:} Let $\displaystyle{ T = \Young(1,2|3,4|5,6)}$. We show the steps to obtain $T^{[2,N]}$ from $T$ using the Sch\"utzenberger algorithm.
$$T =\Young(1,2|3,4|5,6) \quad \rightarrow \quad \Young(,2|3,4|5,6) \quad \rightarrow \quad \Young(2, |3,4|5,6) \quad \rightarrow \quad \Young(2,4|3, |5,6) \quad \rightarrow \quad \Young(2,4|3,6|5)\quad {\rm so }\enspace  T^{[2,N]}=\Young(1,3|2,5|4) $$

\medskip
For $i < j$ we denote by $\pi_{i,j}: \nn_N \rightarrow \nn_{j-i+1}$ the projection which removes all rows and columns numbered $1,2,\ldots,i-1$ or $j+1,\ldots, N$. The image under $\pi_{i,j}$ doesn't depend on the order the rows or columns are removed so it is well-defined (see [M1] 1.3.15). If $\var$ is an orbital variety, the image $\pi_{i,j}(\var)$ is dense in $\overline{\war} \subset \nn_{j-i+1}$ where $\war$ is some orbital variety. We associate to $\war$ the standard Young tableau $\tij$, obtained from $T$ by removing the entries $>j$ and by removing the boxes with entries $<i$ by repeated applications of the above two operations. Again, the order in which these operations are applied doesn't matter, so $\tij$ is well-defined [vanL]. We will use $D_T^{[i,j]}$ to denote the shape of $T^{[i,j]}$.

\bigskip \noindent
{\bf V. The $\sigma$-set of a hypersurface orbital variety}

Let $\var$ be a hypersurface orbital variety and let $\var_R \subset \orb_{\lambda}$ be a Richardson orbital variety such that $\var \subset \overline{\var}_R$ and $\tau(\var) = \tau(\var_R)$. Let $T_R \in \dia_{\lambda}$ be the Richardson tableau (associated to $\var_R$) from which one can obtain a hypersurface tableau $T$ (associated to $\var$). 
We now introduce another subset of $\Pi$, which will be crucial in our study of the non-linear generator $f \in \iv$.

\medskip \noindent
{\bf Definition:} Let $\var$ be a hypersurface orbital variety. The $\sigma$-set $\sigma(\var)$ of $\var$ is the smallest connected subset of $\Pi$ such that the $f$ is contained in $S({\bf n}_{\sigma})$ where $\nn_{\sigma}$ is the subalgebra of $\nn$ generated by $X_{\alpha}$ with $\alpha \in \sigma$. We write $\sigma$ for $\sigma(\var)$ when there is no risk of confusion.

\smallskip
We now study the relationship between $T_R$, $T$ and $\sigma$ by way of the projections $T_R^{[1,N-1]}$, $T_R^{[2,N]}$, $T^{[1,N-1]}$, $T^{[2,N]}$.

\bigskip \noindent
{\bf Theorem 2:} Assume that box $N$ dropped one row to obtain $T$ from $T_R$. Then $\sigma(\var) = \Pi$ if and only if $T_R$ and $\lambda$ satisfy the following two properties:

1. $\mid C_1 \mid $ = $\mid C_{\lambda_1} \mid$. 

2. Let $I:=\mid C_1 \mid $. Then $\lambda_I = \lambda_{I+1}+2$. (In other words, $C_1$ and $C_{\lambda_1}$ are the only chains in $T_R$ of length $I$.)

\medskip \noindent
Proof: Since $\var$ is $\hh$-stable, then the weight of the non-linear defining polynomial $f \in S(\nn)$ is well-defined (every monomial of $f$ has the same weight with respect to $\hh$). By the minimality of $\sigma = \Pi$, this means that for every monomial $m$ of $f$, there exist $i$ and $j$ such that $x_{1i}$ and $x_{jN}$ are factors of $m$.  Since $f$ is the unique nonlinear condition on the $x_{ij}$ ($1\leq i < j\leq N$), the only constraints imposed on the coordinate subsets 
$$\{x_{ij} \enspace \mid  \enspace  1 < i < j \leq N\}
\qquad \hbox{{\rm and}} \qquad \{x_{ij} \enspace \mid  \enspace 1 \leq i < j < N\}$$
(considered as coordinates for either $\var$ or $\var_R$)
are those linear constraints given by: $x_{st} = 0$ if $\alpha(s,t-1) \in R^+(\tau)$. In terms of projections, this translates to $T^{[1,N-1]} = T_R^{[1,N-1]}$ (resp.  $T^{[2,N]} = T_R^{[2,N]}$) and their shapes would be determined uniquely by the restrictions given by $\tau \setminus \{\alpha_1\}$, (resp. $\tau \setminus \{\alpha_{N-1}\}$). 

By definition, $T^{[1,N-1]}$ is obtained by removing box $N$ from $T$. Since $N$ was the block which moved down a row to obtain $T$ from $T_R$, then $T^{[1,N-1]} = T_R^{[1,N-1]}$. On the other hand, $T^{[2,N]} = T_R^{[2,N]}$ means that removing box $1$ from $T$ should precipitate in a shift of boxes which results in box $N$ moving back up one row. When box $1$ is removed from $T$, the boxes corresponding to the rest of the chain $C_1$ move up one row, leaving a space in row $I:= \mid C_1 \mid$ and column $1$. For each remaining box in the $I^{\rm th}$ row, the number to its right is smaller than the number directly underneath it (if it even exists). So, when box $I$ moves up one row, then all remaining boxes in the $I^{\rm th}$ row move over to the left by one space. The remaining boxes in the $l^{\rm th}$ column then move up by one row, where $l$ is the number of boxes in the $I^{\rm th}$ row of $T$. This means that box $N$ is part of this series of shifts if and only if it is in the $l^{\rm th}$ column in $T$. This can happen if and only if $\mid C_1 \mid $ = $\mid C_{\lambda_1} \mid$. Since box $N$ moved down only one row from $T_R$ to $T$, this means that it is in the $I^{\rm th}$ row and $\lambda_I^{\rm th}$ column in $T_R$ and the $(I+1)^{\rm st}$ row and $(\lambda_I - 1)^{\rm st}$ column in $T$ which is true if and only if $\lambda_I = \lambda_{I+1} + 2$.
\qed

From now on, we will call $I$ the {\it thickness} of $\sigma$ and use the notation $I = \thick(\sigma)$.
From Lemma 1 and Theorem 2,
we now know how to ``read off'' $\sigma$ from any tableau $T$ corresponding to a hypersurface orbital variety $\var$. In particular, there exist $i$ and $j$ such that $\pi_{i,j}(\var)$ corresponds to a hypersurface orbital variety in $\nn_{j-i+1}$ whose $\sigma$-set is all of the simple roots for $\nn_{j-i+1}$. The corresponding tableau (obtained by an appropriate renumbering of the entries of $T^{[i,j]}$) will have the same form as that given in Theorem 2. The following corollary also gives the definition for $\thick(\sigma)$ for arbitrary $\sigma \subseteq \Pi$:

\medskip\noindent
{\bf Corollary 1 (Obtaining $\sigma$ and $\thick(\sigma)$ from $T$):} Assume that $T$ was obtained from $T_R$ by dropping the box with the biggest number $j$ of a chain $C_t$ ($1 \leq t \leq \lambda_1$). Let $I :=\mid C_t \mid$. Let $C_s$ ($s<t$) be the chain in $T_R$ of length $I$ such that there is no other length $I$ chain in between $C_s$ and $C_t$. Let $i$ be the smallest number in $C_s$. Then $\sigma(\var) = \{\alpha_i, \alpha_{i+1}, \ldots, \alpha_{j-1}\}$ and $\thick(\sigma):= I= \mid C_s \mid$.

\smallskip\noindent
{\bf Proof:} By Lemma 1, we know that $j$ must be the biggest number in some chain $C_t$ of $T_R$. Let $C_s$ be the previous chain of length $I:=\mid C_t \mid$ in $T_R$ and let $i$ be the smallest number in $C_s$. Then (the renormalized) $C_s$ and $C_t$ correspond to, respectively, the first and last chains of $T_R^{[i,j]}$ and they are the only chains in $T_R^{[i,j]}$ of length $I$. In addition $T^{[i,j]}$ is obtained from $T_R^{[i,j]}$ by dropping the maximal element in the last chain of $T_R^{[i,j]}$. By Theorem 2, $\sigma = \{\alpha_i, \alpha_{i+1}, \ldots, \alpha_{j-1}\}$ and $\thick(\sigma):= I= \mid C_s \mid$.  \qed 

\medskip\noindent
{\bf Example 7:} Consider the hypersurface orbital variety $\var \subset sl(12)$ with tableau
$$T = \Young(1,3,4,7,9,12|2,5,8,10|6|11)$$ We have $\tau(\var) = \{\alpha_1, \alpha_4, \alpha_5, \alpha_7, \alpha_9, \alpha_{10}\}$, so the tableau $T_R$ associated to $\var_R$ is $$T_R = \Young(1,3,4,7,9,12|2,5,8,10|6,11)$$ The chains of $T_R$ are $C_1 = \{1,2\}$, $C_2 = \{3\}$, $C_3=\{4,5,6\}$, $C_4= \{7,8\}$, $C_5=\{9,10,11\}$ and $C_6= \{12\}$. Box $11$ belongs to $C_5$, which has length 3. The preceding chain of length $3$ is $C_3$. So $\sigma(\var) = \{\alpha_4,\alpha_5,\alpha_6,\alpha_7,\alpha_8,\alpha_9,\alpha_{10}\}$ and $\thick(\sigma) = 3$. The hypersurface orbital variety $\war \subset sl(8)$ which is determined by $\tau(\war) = \{ \alpha_1,\alpha_2,\alpha_4,\alpha_6,\alpha_7\}$ (that is $ \sigma \cap \tau(\var) = \{\alpha_4,\alpha_5,\alpha_7,\alpha_9,\alpha_{10}\}$ renormalized under the map $i \rightarrow i-3$) and the same (renormalized) $f$ as $\var$ has standard Young tableau $$T^{[4,11]} = \Young(1,4,6|2,5,7|3|8)$$

\bigskip \noindent
{\bf VI. Van Leeuwen's power-rank conditions}

Now that we have determined the ``location'' of the non-linear generator $f$ from the standard Young tableau for a hypersurface orbital variety, we can work on determining $f$ itself. 

Let $\lambda$ be a partition and let $x \in \orb_{\lambda}$. It is known that $\rk x^j $ is equal to the number of squares beyond the $j^{\rm th}$ column in the Young diagram $D_{\lambda}$. Equivalently, $\dim \ker x^j$ is equal to the number of squares in the first $j$ columns of $D_{\lambda}$. This imposes certain restrictions on the coordinates of $x$. Similar restrictions can be obtained when considering the diagrams associated to the projections $\pi_{i,j}(\nu)$. These restrictions are the so-called {\it power-rank conditions} which were introduced by van Leeuwen [vanL].

Let $\war$ be an orbital variety with tableau $T$ and let $y$ be a 
generic nilpotent matrix in $\war$. For $1 < i\leq j < N$ denote by $y_{[i,j]}:= \pi_{i,j}(y)$ the submatrix of $\nu$ obtained by removing the rows and columns numbered $1,\ldots,i-1$ or $j+1,\ldots,N$. 

\bigskip \noindent
{\bf Theorem 3 [vanL]:} The coordinate vectors of $\nu_{[i,j]}$ satisfy the power-rank conditions imposed by $\dij$. In other words, the coordinate vectors $x_{ij}$ of a generic matrix $\nu \in \war$ satisfies all power-rank conditions imposed by all $\dij$ for $1\leq i\leq j \leq N$.

\bigskip \noindent
{\bf Example 8:} Consider the orbital variety $\var \subset \orb_{(4,2)} \subset \sl(6)$ given by 
$$T = \Young(1,2,4|3,5,6)$$ Notice that $\tau(\var)= \{\alpha_2,\alpha_4\}$ so $T$ is obtained from $T_R$ by dropping box $6$ which is in the last chain. It has length $1$. The only other chain of length $1$ in $T_R$ is $C_1 = \{1\}$ so, in this case, $\sigma = \Pi$ and $I=1$. We show a generic matrix $x \in \nn$ and the associated matrix of $\dij$:

$$x = \left ( \matrix{ 0 & x_{12} & x_{13} & x_{14} & x_{15} & x_{16}\cr
 0 & 0  & x_{23} & x_{24} & x_{25} & x_{26}\cr  0 & 0  & 0 & x_{34} & x_{35} & x_{36}\cr  0 & 0  & 0 & 0 & x_{45} & x_{46}\cr  0 & 0  & 0 & 0 & 0 & x_{56}\cr  0 & 0  & 0 & 0 & 0 & 0 \cr} \right ) \qquad \left 
( \matrix{ {\smallsquares\smallsquares\Young( )} &  {\smallsquares\smallsquares\Young(, )} &  {\smallsquares\smallsquares\Young(, | )} &  {\smallsquares\smallsquares\Young(, , | )} &  {\smallsquares\smallsquares\Young(,, |, )}& {\smallsquares\smallsquares\Young(,, |, , )}\cr
  0  & {\smallsquares\smallsquares\Young( )}  &  {\smallsquares\smallsquares\Young( | )} & {\smallsquares\smallsquares\Young(,| )} & {\smallsquares\smallsquares\Young(,| , )} & {\smallsquares\smallsquares\Young(,,| ,)}\cr  
0  & 0 & {\smallsquares\smallsquares\Young( )}  &  {\smallsquares\smallsquares\Young(, )} &  {\smallsquares\smallsquares\Young(,| )} & {\smallsquares\smallsquares\Young(,,| )} \cr  
0  & 0 & 0 &  {\smallsquares\smallsquares\Young( )} &  {\smallsquares\smallsquares\Young( | )}&{\smallsquares\smallsquares\Young(,| )}\cr  
0  & 0 & 0 & 0 &  {\smallsquares\smallsquares\Young( )}&{\smallsquares\smallsquares\Young( ,)}\cr  0 & 0  & 0 & 0 & 0 & {\smallsquares\smallsquares\Young( )} \cr} \right ) $$
The power-rank condition given by $D_T^{[2,3]} = {\smallsquares\smallsquares\Young(| )}$ says that the matrix $\displaystyle{x_{[2, 3]} = \left ( \matrix{ 0 & x_{23} \cr 0 & 0 \cr} \right )}$ has rank $0$. This implies that $x_{23}=0$ (which we already know since $\alpha_2 \in \tau(\var)$). Likewise, the power-rank condition given by $D_T^{[4,5]}$ forces $x_{45} = 0$. The power-rank condition imposed by $D_T^{[1,6]} = {\smallsquares\smallsquares\Young(,,|, , )}$ says that $\rk x_{[1,6]}^3= 0$. We have that every entry of $x_{[1, 6]}^3$ is $0$ except for the entry in the $1^{\rm st}$ row and $6^{\rm th}$ column which equals $$g:=x_{12} x_{24}x_{46} + x_{12} x_{25}x_{56} + x_{13} x_{34}x_{46} + x_{13}x_{35}x_{56}$$
The condition $\rk x_{[1,6]}^3= 0$ implies that $g=0$. It is easily checked that the remaining power-rank conditions provide {\it trivial power-rank conditions}, that is, they provide no further constraints on the $x_{ij}$. We have $\dim \var = 12 = \dim \var_R - 1$. We conclude that $\var$ is a hypersurface orbital variety and that, since $g$ is irreducible, then $f = g$. Notice that the $\sigma$ that we found in the beginning does correspond to the minimal set such that $f \in S(\nn_{\sigma})$.

\medskip
Let $D_R$ (resp. $D$) denote the shape of $T_R$ (resp. $T$) from the previous example. Then
$$D_R = \Young( , , , | , )\qquad {\rm and} \qquad D = \Young( , , | , , )$$
From $D_R$ and $D$ we have that $\rk x_R^3 = 1$ and $\rk x^3 = 0$. The condition $\rk x^3 = \rk x_R^3 - 1$ results from the fact that a box is dropped down one row in order to obtain $T$ from $T_R$. 

We will use power-rank conditions in precisely the same spirit in order to obtain the non-linear defining generator $f$. Let $\var $ be a hypersurface orbital variety such that $\sigma(\var) = \Pi$ and $I = \thick(\sigma)$. Then $T$ is obtained from $T_R$ by dropping box $N$ from row $I$ to row $I+1$. If box $N$ was in column $k+1$ in $T_R$, then it is in column $k$ in $T$. Let $r$ be the number of boxes after the $k^{\rm th}$ column in $T_R$. Then there are $r-1$ boxes after the $k^{\rm th}$ column in $T$ which simply means that $\rk x^k = \rk x_R^k - 1 = r-1$. Consider the $r \times r$ submatrix $M$ located in the top righthand corner of $x_R^k$. Then $\det M \not= 0$ when $M$ is considered as a submatrix of $x_R^k$, but $\det M = 0$ when it is considered as a submatrix of $x^k$. This means that the non-linear generator $f$ is a factor of $\det M$. In the above example, we had that $f= \det M$. However, this is not always the case as the next example will show:

\medskip\noindent
{\bf Example 9:} Consider the hypersurface orbital variety $\var \subset \orb_{(5,2,2)} \subset sl(9)$ with tableau
$$T = \Young(1,3,6,7,8|2,4|5,9)$$
We have that $\rk x_R^2 = 4$ and $\rk x^2 = 3$. Let $M$ be the $4 \times 4 $ matrix located in the top righthand corner of $x_R^2$. Then
$$\det M = (x_{36}x_{47}-x_{46}x_{37})\cdot g$$
where $g$ is too long to write out. We can check that $g$ is irreducible. Since every monomial in $g$ contains some $x_{1i}$ and some $x_{j9}$ as factors and since $\sigma = \Pi$ in this case, we know that $g$ must equal our nonlinear factor $f \in \iv$.

\medskip \noindent
{\bf Remark:}
Our calculations indicate that $f = \det M$ if and only if the chains in $T_R$ other than the first and last chains all have length less than $I$ {\it or} all have length greater than $I$ (but not both). 

\bigskip \noindent
{\bf VII.  The exact form of the non-linear generator $f \in \iv$}

It is clear that the power-rank conditions do not suffice to give a general formula for $f$. We now present results and a conjecture concerning its exact form. Without loss of generality, we can restrict ourselves to the case $\sigma = \Pi$. 
 A generic matrix $x_R \in \var_R$ has the 
following form:

$$x_R = \left ( \matrix{0 & \cdots& 0 & x_{1,I+1} & \cdots& \cdots & \cdots & \cdots &
x_{1,N}\cr
:& & :&:& :& : & : & : & :\cr
: & & 0 & x_{I,I+1} & \cdots & \cdots  & x_{I,N-I+1} &\cdots & :\cr
: & &  &0 & x_{I+1,I+2} & & : & & :\cr
: & &  & 0 & & :& & &:\cr
: & & & & 0 & & : & &: \cr
: & & & & &x_{N-I-1,N-I} & : & & : \cr
0&  & & \cdots &  & 0 & x_{N-I,N-I+1} & \cdots & x_{N-I,N}\cr 
0 & & & & \cdots& 0& 0 & \cdots & 0\cr
: & & & & & & & & :\cr
0& \cdots & \cdots & \cdots & \cdots & \cdots & \cdots & \cdots & 0\cr} \right ).$$
where $x_{k,l+1}=0$ if and only if $\alpha(k,l) \in R^+(\tau)$. 

For a matrix $y$, denote by $\cMin_k(y)$ the $(N-k) \times (N-k)$ submatrix in the top right corner of $y$.
Let $t \in \CC$. Then
$$\cMin_{I}(x_R + t\id) = \left ( \matrix{x_{1, I+1} & \cdots& \cdots & \cdots & \cdots & &
x_{1,N}\cr
:& :& : & : & : & & :\cr
x_{I, I+1} & \cdots & \cdots  & \cdots &\cdots & & :\cr
t & x_{I+1, I+2} & & : & & :\cr
0 & t & & : & &:\cr
: & 0 & & : & &: \cr
: & &t &x_{N-I-1,N-I} & : & & : \cr
0& \cdots & 0 & t & x_{N-I, N-I+1} & \cdots & x_{N-I ,N}\cr } \right )$$
We have
$$\det(\cMin_{I}(x_R + t\id)) =  m_{N-I} + m_{N-I-1}t + \cdots + m_It^{N-2I}$$
where the $m_j$ are (up to sign) sums of $j\times j$ minors in $x$. 
Let $k \geq 0$ and let $\ui := (i_1,i_2,\ldots,i_k)\in \NN^k$ satisfy $I+1 \leq i_1 < i_2 < \cdots < i_k \leq N-I$. Consider the $(k+I) \times (k+I)$ matrix
$$A_{\ui} = \left ( \matrix{
x_{1,i_1} & x_{1,i_2} & \cdots &\cdots &\cdots& x_{1,N-I+1} & \cdots & x_{1,N}\cr
\vdots & \vdots & & & & & \vdots \cr
x_{I,i_1} & x_{I,i_2} & \cdots &\cdots &\cdots& x_{I,N-I+1} & \cdots & x_{I,N}\cr
0 & x_{i_1,i_2} & x_{i_1,i_3} & \cdots & x_{i_1,i_k} & x_{i_1,N-I+1}& \cdots & x_{i_1,N}\cr
0 & 0 & x_{i_2,i_3} & & & & & \vdots \cr
\vdots & & & & & & & \vdots \cr
0 & \cdots & \cdots & \cdots & \cdots & x_{i_k,N-I+1} & \cdots & x_{i_k,N} \cr} \right ) $$
Then $m_{k+I} = \sum_{\ui} \det(A_{\ui})$ where the sum is over all possible $k$-tuples $\ui$.

\bigskip\noindent
{\bf Lemma 2:} $m_{k+I} = 0$ if and only if $k+I > \lambda_1 + \cdots + \lambda_I - I$.

\noindent 
{\bf Proof:} Any two $k$-tuples $\ui$, $\uj$ differ by at least one entry. Hence any monomial in $\det A_{\ui}$ differs from any monomial in $\det A_{\uj}$ by at least one factor $x_{st}$. So the algebraic independence of the $x_{ij}$ implies that there can be no cancelations of monomials from determinants of different $A_{\ui}$. Hence $m_{k+I}$ is non-zero if and only if there exists an $A_{\ui}$ such that $\det(A_{\ui}) \not= 0$.

Notice that $x_{ij} = 0$ if and only if $x_{rs} = 0 $ for all $(r,s)$ such that $r \geq i$ and $s \leq j$ (all coordinates below or to the left of $x_{ij}$). If at least one of the diagonal elements of $A_{\ui}$ is zero, then $A_{\ui}$ is an upper-triangular block matrix where at least one of the block matrices has a zero column or zero row. In this case, $\det(A_{\ui})=0$. Clearly, when all of its diagonal elements are non-zero, we have that $\det(A_{\ui}) \not=0$. Therefore, $\det(A_{\ui}) \not=0$ if and only if all of its diagonal coordinates are non-zero. When $k \leq I $, all diagonal elements are of the form $x_{ij}$ where $1 \leq i \leq I$ or $N-I+1 \leq j \leq N$ so they are all non-zero and $m_{k+I} \not=0$. 

Now consider $k \geq I+1$. Then $A_{\ui}$ has $k-I$ diagonal elements of the form $x_{i_j,i_{j+I}}$. For any given $k$-tuple $\uj$, let $d_{\uj}$ denote the number of non-zero diagonal elements of the form $x_{j_t, j_{t+I}}$ where $1 \leq t \leq k$. Let $c_k := \max d_{\uj}$ where the maximum is taken over all $k$-tuples $\uj$. Let $c:= \max_{k} c_k$ where the maximum is taken over all $k$. 
We have that $\det(A_{\ui}) =0$ for all $k$-tuples $\ui$ if and only if 
$k-I > c$. This maximum $c$ is attained at the -tuple which contains every non-zero coordinate of the form $x_{i,i+I}$. In other words, $c$ equals the number of positive roots of length $I$ in $R^+(\omega) \setminus R^+(\tau\cap\omega)$ where $\omega :=\{\alpha_{I+1}, \alpha_{I+2}, \ldots, \alpha_{N-I-1}\}$. There are $N-3I$ positive roots of length $I$ in $R^+(\omega)$. The number of positive roots of length $I$ in $R^+(\tau\cap\omega)$ is the number of boxes after the $I^{\rm th}$ row of $T_R$. In other words, it is $N- (\lambda_1 + \cdots + \lambda_I)$. Therefore, $c = \lambda_1 + \lambda_2 + \cdots + \lambda_I - 3I$. We therefore have that $ m_{k+I} = 0$ if and only if $k+I > \lambda_1 + \cdots + \lambda_I -I$.  \qed

\noindent
For ease of notation, set $l(\lambda) := \lambda_1 + \cdots + \lambda_I - I$. For any $\hh$-semiinvariant polynomial $p \in S(\nn)$, we denote the weight of $p$ by $\wt(p)$. Let $x_{l_1 l_2}x_{l_3 l_4}\cdots x_{l_{s-1}l_s}$ be any monomial term of $p$. Then $\wt(p) = \wt(x_{l_1 l_2}x_{l_3 l_4}\cdots x_{l_{s-1}l_s}) = \wt(x_{l_1 l_2}) + \wt(x_{l_3 l_4})+\cdots +\wt(x_{l_{s-1}l_s})$. We have that $\wt(x_{ij}) = \alpha(i,j-1)$ for $1\leq i< j \leq N$.

\bigskip\noindent
{\bf Proposition 1:} When $I=1$ then $f = m_{l(\lambda)} = m_{\lambda_1 - 1}$.

\noindent
{\bf Proof:} In this case, the matrices $A_{\ui}$ corresponding to the $(\lambda_1 - 2)$-tuples $\ui$ are diagonal and 
$$\det(A_{\ui}) = x_{1,i_1}x_{i_1,i_2}\cdots x_{i_{\lambda_1-2},N}$$
It follows that $m_{\lambda_1 - 1}$ is the sum of all possible monomials of this type. Now, consider 
$$x_R = \left ( \matrix{0 & A & B \cr 0 & C & D \cr 0& 0 & 0\cr} \right ) \in  \var_R $$
where $A$ is the $1 \times (N-2)$ row vector $[x_{1,2} \enspace \cdots \enspace x_{1,N-1}]$, $B$ is the $1 \times 1$ matrix $[x_{1,N}]$,  $D$ is the $(N-2) \times 1$ column vector $[x_{2,N}\enspace  \cdots \enspace  x_{N-1,N}]^T$, $C$ is the $(N-2) \times (N-2)$ square matrix that is left. The $0$s represent zero matrices of the appropriate size. For all $j \geq 2$, we have
$$
x^{j}_R = \left ( \matrix{0 & AC^{j-1} & AC^{j-2}D \cr 0 & C^{j} & C^{j-1}D \cr 0 & 0 & 0\cr} \right ). $$
Notice that $C = \pi_{[2,N-1]}(x_R)$. In other words, $C$ is the generic matrix associated to the Richardson orbital variety with standard Young tableau $T_R^{[2,N-1]}$. The tableau $T_R^{[2,N-1]}$ is obtained from $T_R$ by removing the first and last chains, which correspond simply to box $1$, resp. box $N$. This means that there are only $\lambda_1 - 2$ columns in $T_R^{[2,N-1]}$. By the power-rank conditions, $\rk(C^{\lambda_1 - 3}) > 0$, but $\rk(C^{\lambda_1 - 2}) = 0$ and $\rk(C^{\lambda_1 - 1}) = 0$. So $C^{\lambda_1 - 3} \not= 0$, but $C^{\lambda_1 - 2} = C^{\lambda_1 - 1} = 0$.

There is only one box (which is box $N$, in fact) past the $(\lambda_1 -1)^{\rm st}$ column in $T_R$. Therefore, when $j = \lambda_1 - 1$, the matrix $x^j_R$ has rank $1$. Since $C^{\lambda_1 - 1} = 0$ and $C^{\lambda_1 -2} = 0$, we have that $x^{\lambda_1 -1}_R$  is zero everywhere except at the $1 \times 1$ matrix $A C^{j-2}D$ in the top right corner. Now, $\rk(x^{\lambda_1-1}) = \rk(x^{\lambda_1 -1}_R) -1 = 0$ so $A C^{j-2}D=0$ when considered as an entry in $x^{\lambda_1 -1}$. Hence $f$ divides $ A C^{j-2} D$. This implies that $\wt( A C^{j-2} D) - \wt(f) $ is a positive sum of positive roots.

We have $\wt(A C^{j-2} D) = \wt(x_{1,i_1}x_{i_1,i_2}\cdots x_{i_{\lambda_1-2},N}) = \alpha(1,N-1)$. 
We will show that $f= A C^{j-2} D$ by showing that $\wt(f)=\wt( A C^{j-2} D) = \alpha(1,N-1)$. By the definition of $\sigma$, we know that the coefficient of $\alpha_1$ in $\wt(f)$ is non-zero. If $\wt(f) \not= \alpha(1,N-1)$ then there is a smallest $i$ ($1 < i < N-1$) such that the coefficient of $\alpha_i$ in $\wt(f)$ is $0$. We can then write $\wt(f) = \alpha(1,i-1) + \beta$ where $\alpha_i$ is not a summand of $\beta$. Since every chain (besides the first and last) in $T_R$ is of length $\geq 2$, then $\alpha_{i-1} \in \tau$ or $\alpha_i \in \tau$ (or both). But, 
$$s_{i-1}(\wt(f)) = s_{i-1}(\alpha(1,i-1) + \beta) = \alpha(1,i-2) + \beta \not= \wt(f).$$
Similarly, $s_i(\wt(f)) = \alpha(1,i) + s_i(\beta)\not= \wt(f)$ (independently of whether or not $\alpha_{i+1}$ is a summand of $\beta$). Therefore, $\wt(f)$ is not invariant under the action of $\tau$, which brings us to a contradiction. Finally, since $A C^{j-2} D = \sum_{\ui}x_{1,i_1}x_{i_1,i_2}\cdots x_{i_{\lambda_1-2},N} = m_{\lambda_1 - 1}$ we have $f = m_{\lambda_1 - 1}$.
 \qed

\medskip\noindent
{\bf Example 10:} In the case of Example 8, we have that
$$\cMin_1(x_R + t\id) =  \left ( \matrix{x_{1 2} & x_{1 3}& x_{1 4} & x_{15} & x_{16}\cr
t& 0 & x_{2 4} & x_{2 5} & x_{26}\cr
0 & t & x_{34} & x_{35} &x_{36}\cr
0 & 0 & t & 0& x_{46}\cr
0 & 0 &0 & t & x_{56}\cr} \right )$$
We have $\det \cMin_1(x_R + t\id) =x_{16} t^4 - (x_{12}x_{26}+ x_{13}x_{36} + x_{14}x_{46} + x_{15}x_{56})t^3 + (x_{12} x_{24}x_{46} + x_{12} x_{25}x_{56} + x_{13} x_{34}x_{46} + x_{13}x_{35}x_{56}) t^2$
and indeed $m_3 = x_{12} x_{24}x_{46} + x_{12} x_{25}x_{56} + x_{13} x_{34}x_{46} + x_{13}x_{35}x_{56}$ is the generator that we had previously found.

\bigskip\noindent
{\bf Proposition 2:} If $\lambda_{I+1} = 0$ then $f = m_{l(\lambda)}$.

\noindent
{\bf Proof:} For ease of notation, we will denote the coordinates $x_{i,j+I}$ by $M_{i,j}$ for $1\leq i,j \leq N-I$. 

Since $\lambda_{I+1} = 0$, then $\lambda_{I} = 2$. Consequently, the only chains of length $\geq I$ in $T_R$ are the first and last chains. This means that the subdiagonal coordinates $M_{i, i-1}$ are not identically $0$. We show that this implies that $d:= \det (\cMin_I(x_R))$ is irreducible. In fact, suppose that $d = pq$ where both $p$ and $q \in \CC[M_{i,j}\mid 1 \leq i, j \leq N-I]$. We claim that $p$ and $q$ are functions on disjoint sets of row vectors. In fact, suppose that they aren't. Then, for some row $i$ and some columns $j$ and $k$, we have that both $p$ and $q$ are functions of the coordinates $M_{i,j}$ and $M_{i,k}$. This means that we can write $p=M_{i,j}p_1 + M_{i,k}p_2 + p_3$ and $q=M_{i,j}q_1 + M_{i,k}q_2 + q_3$ where the $p_l$ and $q_l$ ($l = 1,2,3$) are in $\CC[M_{i,j}\mid 1 \leq i, j \leq N-I]$ and where $p_1$ and $q_1$ do not depend on $M_{i,j}$, where $p_2$ and $q_2$ do not depend on $M_{i,k}$ and where $p_3$ and $q_3$ do not depend on either $M_{i,j}$ or $M_{i,k}$. Then $pq = M_{i,j}^2p_1q_1 + M_{i,k}^2p_2q_2 + M_{i,j}M_{i,k}(p_1q_2+p_2q_1) + r$ where the degree of either $M_{i,j}$ or $M_{i,k}$ in any term of $r$ is at most $1$. The degrees of $M_{i,j}$ and $M_{i,k}$ are at most $1$ in every term of $d$ and no term of $d$ contains $M_{i,j}M_{i,k}$. Therefore we have that $p_1q_1= p_2q_2=p_1q_2+p_2q_1=0$. Since $ \CC[M_{i,j}\mid 1 \leq i, j \leq N-I]$ is a domain, then either $p_1 = 0 $ or $q_1=0$. Without loss of generality, we can assume that $p_1=0$. Then either $p_2=0$ or both $q_1=q_2=0$. In either case, one of the factors $p$ or $q$ does not depend on $M_{i,j}$ and $M_{i,k}$. So $p$ and $q$ depend on disjoint sets of row vectors. A similar argument shows that $p$ and $q$ depend on disjoint sets of column vectors. Therefore, there exist subsets $G$ and $H$ in $\{1,2,\ldots,N-I\}$ such that $p \in \CC[M_{i,j}\mid i\in G, j\in H]$ and $q \in \CC[M_{i,j} \mid i \in G^{\prime}, j \in H^{\prime}]$ where $G^{\prime}$ and $H^{\prime}$ are the complements of $G$ and $H$ in $\{1,2,\ldots,N-I\}$. Without loss of generality, we can assume that $G$ is non-empty. We claim that $G=H$. In fact, consider the following specialization $M$ of $\cMin_I(x_R)$: set $M_{i,j}=0$ for $i\not=j$. Then $d(M) = M_{1,1}M_{2,2}\cdots M_{N-I ,N-I} = p(M)q(M)$. This means that $i \in G$ if and only if $i \in H$. So $G=H$. Now we claim that $G= \{1,2,\ldots, N-I\}$. In fact, consider the following specialization $P$ of $\cMin_I(x_R)$: set $M_{i,j}=0$ for $(i,j)$ not of the form $(i,i+1)$ or $(1,N-I)$. Then $d(P) = M_{1,2}M_{2,3}\cdots M_{N-I-1, N-I}M_{1, N-I} = p(P)q(P)$. This means that if some $i \in G$ then $i+1 \in H=G$. Therefore, $G=H=\{1,\ldots, N-I\}$ and $G^{\prime}=H^{\prime}=\emptyset$ which means that $q$ is constant. Therefore, $d = \det (\cMin_I(x_R))$ is irreducible.

The tableau $T$ is obtained from $T_R$ by dropping box $N$ from the $I^{\rm th}$ row and $2^{\rm nd}$ column in $T_R$ to the $(I+1)^{st}$ row and $1^{st}$ column in $T$. We have that $\rk(x) = \rk(x_R)-1$ and, by irreducibility of $\det(\cMin_I(x_R))$, we have $f= \det(\cMin_I(x_R))$. On the other hand, there are no boxes in $T_R$ past the $I$th row, so $ {l(\lambda)} = rk(x_R) = N-I$. Therefore, $m_{l(\lambda)} = \det(\cMin_I(x_R))$. \qed

\medskip\noindent
{\bf Example 11:} Consider the hypersurface orbital variety $\var \subset \orb_{(3,1,1)} \subset sl(5)$ with tableau
$$T = \Young(1,3,4|2|5)$$
We have $\sigma = \Pi$, $I=2$ and
$$\cMin_2(x_R + t\id) = \left ( \matrix{x_{13} & x_{14}& x_{15}\cr 
x_{23} & x_{24} & x_{25}\cr
t & x_{34} & x_{35}\cr } \right ). $$
Then $\det \cMin_2(x_R + t\id) = m_2 t + m_3$. We have that $m_3$ is irreducible. In addition $m_3 = 0$ is exactly that constraint given by $\rk x = 2$. So $m_3$ is the nonlinear generator in $\iv$.

\medskip
The results of the previous two propositions and of explicit calculations using MAPLE for most cases up to $\SL(14)$ have led us to believe that $f$ is {\it always} equal to $m_{l(\lambda)}$. We claim:

\medskip\noindent
{\bf Conjecture:} Let $\var$ be a hypersurface orbital variety. Let $\var_R$ be the Richardson orbital variety with the same $\tau$-invariant as $\var$. Assume that $\sigma(\var) = \{\alpha_i, \ldots, \alpha_j\}$ and let $I:= \thick(\sigma(\var))$. Let $x_R$ represent a generic matrix in $\var_R$. Let $m$ be the coefficient of the smallest power of $t$ in $\det \cMin_I((x_R + t\id)_{[i,j]})$. Then $f=m$ and $\iv = \langle X_{\alpha}, \enspace \alpha \in R^+(\tau(\var)),\enspace m \rangle$.

\medskip\noindent
{\bf Example 12:} Let $\var$ be as in Example 7. Then $\sigma = \{\alpha_4,\alpha_5,\alpha_6,\alpha_7,\alpha_8,\alpha_9,\alpha_{10}\}$ and $\thick(\sigma) = 3$. In this case 
$$\cMin_3((x_R+t\id)_{[4,10]}) = \left ( \matrix{x_{4 7} & x_{4 8} & x_{4 9} & x_{4 \thinspace 10} & x_{4 \thinspace 11} \cr x_{5 7} & x_{5 8} & x_{5 9} & x_{5 \thinspace10} & x_{5\thinspace 11} \cr
 x_{6 7} & x_{6 8} & x_{6 9} & x_{6 \thinspace10} & x_{6 \thinspace11} \cr
 t & 0 & x_{7 9} & x_{7 \thinspace10} & x_{7 \thinspace11} \cr
0 & t & x_{8 9} & x_{8 \thinspace10} & x_{8 \thinspace11} \cr} \right )$$
and $\det \cMin_3((x_R+t\id)_{[4,10]}) = m_3 t^2 + m_4 t + m_5$. We have $f = m_5$. So $\iv = \langle f, \thinspace x_{ij} \thinspace \mid \thinspace \alpha(i,j-1) \in R^+(\tau) \rangle$.

\bigskip\noindent
{\bf VIII. The characteristic polynomial $p_{\var}$}

\medskip
To every orbital variety $\var$, one can associate its {\it characteristic polynomial} $p_{\var} \in S(\hh^*)$ ([J1], [BBM]). This $W$-harmonic polynomial has degree = $\vert R^+ \vert - \dim \var$ [J2]. Although $p_{\var}$ is explicitly known
for a number of examples, there is no known (explicit) general
formula. However, when $\var$ is a hypersurface orbital variety, we can use the following theorem to give an explicit formula for $p_{\var}$. This theorem is an immediate consequence of work by Joseph (see [J1, 82.], [J2, 2.9]). A different proof is given by Borho, Brylinski and MacPherson [BBM, 4.15]. 

\medskip\noindent
{\bf Theorem 4:} Let $\var$ be a complete intersection of codimension $d$ in $\nn$, defined by homogeneous equations $f_1$, $f_2, \ldots, f_d \in S(\nn)$ of weights $\mu_1, \dots, \mu_d$. Then $p_{\var} = \mu_1 \mu_2 \cdots \mu_d$.

\medskip
In our case, the formula for $\wt(f)$, where $f=0$ is the unique non-linear defining equation of $\var$, is an easy consequence of Propositions 1 and 2 (or the Conjecture, once proven). Since $\det(\cMin_I(x_R+t\id))$ is homogeneous and $\wt(t) = 0$ we have that $\wt(m_I) = \wt(m_{I+1}) = \cdots = \wt(m_{l(\lambda)}) = \wt(f)$. $m_I$ is simply the determinant of the top right $I \times I$ minor in $x_R$ and its weight is equal to the weight of any of its summands.  Therefore, 
$$\eqalign{\wt(f) = \wt(m_I) &= \wt(x_{1,N}x_{2,N-1}\cdots x_{I,N-I+1})\cr
&=\wt(x_{1,N}) + \wt(x_{2,N-1})\cdots + \wt(x_{I,N-I+1})\cr
&=\alpha(1,N-1)+\alpha(2,N-2) + \cdots + \alpha(I,N-I)\cr}$$ and we get the desired formula. Clearly this formula holds for $\var$ satisfying the assumptions of Propositions 1 and 2.

\medskip\noindent
{\bf Corollary 2:} Let $\var$ be a hypersurface orbital variety with $\tau$-invariant $\tau(\var)=\tau$, $\sigma$-set $\sigma = \{\alpha_i, \alpha_{i+1}, \ldots, \alpha_{j-1}\}$ and $I = \thick(\sigma)$. Then $$p_{\var} = (\prod_{\alpha \in R^+(\tau)} \alpha) \cdot (\sum_{k=0}^{I-1} \alpha(i+k,j-k))$$

\medskip\noindent
{\bf Example 13:} Let $\var$ be as in Examples 7 and 12. Then $$\eqalign{\wt(f) = \wt(x_{69}) + \wt(x_{5 \thinspace 10}) + \wt(x_{4 \thinspace 11}) &= \alpha(6,8) + \alpha(5,9) + \alpha(4,10)\cr
&= \alpha_4 + 2 \alpha_5 + 3\alpha_6 + 3\alpha_7 + 3\alpha_8 + 2\alpha_9 + \alpha_{10}\cr}$$
Therefore, $p_{\var} = \alpha_1 \alpha_4 \alpha_5 \alpha_7 \alpha_9 \alpha_{10}(\alpha_4 + \alpha_5)(\alpha_9 + \alpha_{10}) (\alpha_4 + 2 \alpha_5 + 3\alpha_6 + 3\alpha_7 + 3\alpha_8 + 2\alpha_9 + \alpha_{10})$.

\bigskip\noindent
{\bf References}

[Be] E. Benlolo, Sur la quantification de certaines vari\'et\'es orbitales. {\it Bull. Sci. Math.} {\bf 118} (1994), no. 3, 225-243. 

[BV] D. Barbasch and D. Vogan, Primitive ideals and orbital
  integrals in complex classical groups, {\it Math. Ann.} {\bf 259}
  (1982) 153-199.

[B] W. Borho, Nilpotent orbits, primitive ideals, and
  characteristic classes (a survey), {\it Proceedings of the
  International Congress of Mathematicians, Vol. 1,2 (Berkeley,
  Calif., 1986)}, Amer. Math. Soc, Providence, (1987) 350-359.

[BBM] W. Borho, J-L Brylinski and R. MacPherson, Equivariant K-Theory Approach to Nilpotent Orbits, IHES/M/86/13

[He] W. Hesselink, Singularities in the nilpotent scheme of a classical group, {\it Trans. Am. Math. Soc.} {\bf 222} (1976) 1 - 32.

[Ho] A. Hotta. On Joseph's construction of Weyl group
representations. {\it Tohuku Math. J.} {\bf 36} (1984), 49--74.

[Ja] J. C. Jantzen, {\it Einh\"ullende Algebren halbeinfacher Lie-Algebren}, Ergebnisse der Mathematik und ihrer Grenzgebiete, {\bf 3}, Springer-Verlag, 1983.

[J1] A. Joseph, On the characteristic polynomials of orbital varieties, {\it Ann. scient. Ec. Norm. Sup}, 4eme s\'erie, t.22 (1989), 569 - 603.

[J2] A. Joseph, On the Variety of a Highest Weight Module, {\it J. Algebra}, {\bf 88}, No. 1 (1984), 238-278.

[J3] A. Joseph, Enveloping Algebras: Problems Old and
New, 385-413, {\it Progress in Mathematics}, {\bf 123}, Birkh{\"a}user,
Boston, 1994.

[J4] A. Joseph, Orbital varietes of the minimal orbit. {\it Ann. Sci. \'Ecole Norm. Sup. (4)} {\bf 31} (1998), no. 1, 17-45. 

[Mc] W. M. McGovern, Dixmier algebras and the orbit
  method. {\it Operator algebras, unitary representations, enveloping
  algebras, and invariant theory (Paris, 1989)}. Progr. Math, 92,
  Birkh{\"a}user, Boston, 1990.

[M1] A. Melnikov, Orbital Varieties and Order Relations on Young Tableaux, (1995) preprint.

[M2] A. Melnikov, Orbital Varieties in $\sl(n)$ and the Smith Conjecture, {\it J. Algebra} {\bf 200} (1998)  1 - 31.

[M3] A. Melnikov, Irreducibility of the associated varieties of simple highest weight modules in $sl(n)$. {\it C. R. Acad. Sci. Paris Sér. I Math.} {\bf 316} (1993), no. 1, 53-57. 

[Sp] N. Spaltenstein, Classes unipotentes de sous-groupes de Borel. Lecture Notes in Mathematics, {\bf 964} Springer-Verlag, Berlin-New York, 1982.

[SS] T. A. Springer, R. Steinberg, {\it Conjugacy classes}, Lecture Notes in Mathematics, {\bf 131}, Springer, 1970, 167-266.

[St1] R. Steinberg, An Occurrence of the Robinson-Schensted Correspondence, {\it J. Algebra}, {\bf 113}, (1988), 523-528.

[St2] R. Steinberg, On the desingularization of the unipotent variety, {\it Invent. Math.}, {\bf 36} (1976), 209-224.

[vanL] M. A. A. van Leeuwen, The Robinson-Schensted and Sch\"utzenberger algorithms, Part II: Geometric interpretations, CWI report AM-R9209 (1992). This report is also available electronically via {\bf http://www.cwi.nl/cwi/publications/\#AM}.

\bye